\documentclass[12pt]{amsart}
\usepackage{amsmath,amssymb,graphicx, amscd}
\usepackage[all]{xy}
\usepackage[bbgreekl]{mathbbol}
\usepackage{verbatim}
\usepackage{enumerate}

\newtheorem{thm}{Theorem}[section]

\newtheorem{lem}[thm]{Lemma}
\newtheorem{ass}[thm]{Assumption}
\newtheorem{prop}[thm]{Proposition}
\theoremstyle{definition}
\newtheorem{defn}[thm]{Definition}
\theoremstyle{remark}

\newtheorem{rem}[thm]{Remark}
\numberwithin{equation}{section}

\DeclareFontFamily{U}{rsf}{} \DeclareFontShape{U}{rsf}{m}{n}{
  <5> <6> rsfs5 <7> <8> <9> rsfs7 <10->  rsfs10}{}
\DeclareMathAlphabet{\mathscr}{U}{rsf}{m}{n}

\newcommand{\abs}[1]{\left\vert#1\right\vert}

\newcommand{\A}{\mathcal{A}}

\newcommand{\Cc}{C}
\newcommand{\Dc}{D}

\newcommand{\Pc}{{\mathcal P}}

\newcommand{\Zb}{{\mathbb Z}}

\newcommand{\As}{\text{\bfseries\sf{A}}}
\newcommand{\Rs}{\text{\bfseries\sf{R}}}

\renewcommand{\imath}{\sqrt{-1}}

\newcommand{\Eb}{\mathbb E}

\newcommand{\Mb}{\mathbb M}

\newcommand{\Nb}{\mathbb N}

\newcommand{\Xb}{\mathbb X}

\newcommand{\np}{ }
\newcommand{\dbar}{\overline{\partial}}
\DeclareMathOperator{\modl}{Mod} 
 
\DeclareMathOperator{\Hom}{Hom}
\DeclareMathOperator{\Ext}{Ext} 

\DeclareMathOperator{\Ob}{Ob} 
\DeclareMathOperator{\Ho}{Ho}

\newcommand{\mA}{\mathcal{A}}
\newcommand{\mB}{\mathcal{B}}
\newcommand{\mC}{\mathcal{C}}
\newcommand{\mD}{\mathcal{D}}
\newcommand{\mX}{\mathcal{X}}

\newcommand{\tA}{\text{\bfseries\sf{A}}}
\newcommand{\tB}{\text{\bfseries\sf{B}}}
\newcommand{\tC}{\text{\bfseries\sf{C}}}
\newcommand{\tD}{\text{\bfseries\sf{D}}}

\newcommand{\tf}{\text{\bfseries\sf{f}}}
\newcommand{\tg}{\text{\bfseries\sf{g}}}
\newcommand{\tk}{\text{\bfseries\sf{k}}}
\newcommand{\tl}{\text{\bfseries\sf{l}}}

\newcommand{\tE}{\text{\bfseries\sf{E}}}
\newcommand{\tF}{\text{\bfseries\sf{F}}}

\newcommand{\tS}{\text{\bfseries\sf{S}}}
\newcommand{\tM}{\text{\bfseries\sf{M}}}
\newcommand{\tN}{\text{\bfseries\sf{N}}}
\newcommand{\tH}{\text{\bfseries\sf{H}}}
\newcommand{\tX}{\text{\bfseries\sf{X}}}

\newcommand{\mM}{\mathcal{M}}
\newcommand{\mN}{\mathcal{N}}
\newcommand{\mE}{\mathcal{E}}
\newcommand{\mF}{\mathcal{F}}
\newcommand{\mH}{\mathcal{H}}
\newcommand{\mR}{\mathcal{R}}
\newcommand{\mP}{\mathcal{P}}
\newcommand{\mQ}{\mathcal{Q}}
\newcommand{\mG}{\mathcal{G}}

\numberwithin{equation}{section}

\begin{document}

\title[Milnor Descent For Cohesive dg-categories]{Milnor Descent For Cohesive dg-categories}
\author{Oren Ben-Bassat, Jonathan Block}\thanks{J.B. partially supported by NSF grant DMS10-07113}%
\address{Oren Ben-Bassat, Department of Mathematics, University of Haifa, Haifa, Israel.\newline \indent Jonathan Block, Department of Mathematics, University of Pennsylvania, Philadelphia, PA, USA}
\email{ben-bassat@math.haifa.ac.il, blockj@math.upenn.edu}%
\dedicatory{}
\subjclass{46L87, 58B34}%
\thanks{}
\keywords{}%

\begin{abstract}

\noindent  We show that the functor from curved differential graded algebras to differential graded categories, defined by the second author in \cite{B}, sends Cartesian diagrams to homotopy Cartesian diagrams, under certain reasonable hypotheses.   This is an extension to the arena of dg-categories of a construction of projective modules due to Milnor.  As an example, we show that the functor satisfies descent for certain partitions of a complex manifold. 
\end{abstract}
\maketitle 
\section{Introduction}  
Differential graded categories \cite{Ke2} (and their relatives $A_\infty$ and $\infty$-categories) are of increasing interest and value in many areas of algebraic, symplectic, differential and noncommutative geometry. They arise as better behaved enhancements of derived (or triangulated) categories. They provide a context in which  to understand various relations amongst disparate objects. A classic example along this line is the Riemann-Hilbert correspondence between $D$-modules and perverse sheaves. A more recent one is Mirror Symmetry, a symmetry coming out of physics that relates the symplectic geometry of one Calabi-Yau manifold to the complex geometry of a mirror Calabi-Yau manifold. By now the examples are too numerous to survey completely.

In \cite{B} the second author introduced a functor from curved differential graded algebras to differential graded categories
\[
\As\mapsto \Pc_\As .
\] This provided a way to ``present" many dg-enhancements of well-known categories of interest in geometry as well many new ones. Further, the beginnings of a calculus was developed to be able to calculate with $\Pc_\As$. In this paper we continue the development of this calculus by proving a descent theorem.  Of course, descent was one of the original motivations of considering dg-categories and descent results have been proven for dg-enhancements of the category of perfect complexes of sheaves on a scheme $X$, such as what is called $L_{pe}(X)$ in \cite{T}. 

\begin{prop}(A. Hirschowitz, C. Simpson, see \cite{T}, Proposition 11)
 Let $X = U\cup  V$, where $U$ and $V$ are two Zariski open subschemes. 
Then the following square 
\begin{equation}
\begin{array}{ccc}
L_{pe} (X ) & \to & L_{pe} (U )\\ 
\downarrow & & \downarrow \\
L_{pe} (V ) & \to & L_{pe} (U \cap V ) 
\end{array}
\end{equation}
is homotopy cartesian in the model category  dg-cat  of dg-categories. 
\end{prop}

Our descent principle  will imply a similar result for certain covers of a complex manifold $X$ which can look very different than a cover by open sets.  From an algebraic perspective, our theorem is an extension to the world of dg-categories of a theorem (see Subsection \ref{CatProjMod}) of Milnor \cite{M} which recovers a projective module over the fiber product of rings $R \times_{T} S$  from a pair of projective modules on $R$ and $S$ and a gluing isomorphism over $T$, under the condition that one of the maps to $T$ is surjective. Our main result (Theorem \ref{main}) takes the form of a homotopy fiber product of dg-categories 
\begin{equation}
\begin{array}{ccc}
\mathcal{P}_{\tA} & \to &\mathcal{P}_{\tB}\\ 
\downarrow & & \downarrow \\
\mathcal{P}_{\tC} & \to &\mathcal{P}_{\tD}
\end{array}
\end{equation}
associated to a diagram of dg-algebras 
\begin{equation}
\begin{array}{ccc}
\tA & \to &\tB \\ 
\downarrow & & \downarrow \\
\tC & \to & \tD
\end{array}
\end{equation}
satisfying certain assumptions. One also take curved dg-algebras instead of dg-algebras. We outline some geometric applications of this theorem to complex manifolds in Section \ref{CCM}. In these examples, the dg-algebras involved are the Dolbeaut dg-algebras of the complex manifolds (or direct limits of such). We should emphasise that our theorem involves descent for certain dg-categories associated to these dg-algebras and not descent for the algebras themselves. In Section \ref{HFP} we give some Mayer-Vietoris sequences for some different type of invariants of dg-categories appearing in homotopy fiber product diagrams. Since the first version of this article as put on the arxiv \cite{BB}, some of the methods presented in it have found applications to the study of formal neighborhoods \cite{Yu}. See also \cite{BT} for a similar type of descent in the algebraic context.

\section{Review of Basic Definitions}
In this section we recall the basic definitions and theorems from \cite{B}.

\begin{defn} For complete definitions and facts regarding dg-categories, see \cite{BK}, \cite{Dr}, \cite{Ke} and \cite{Ke2}. Fix a field $k$. A {\em differential graded category (dg-category) } is a
category enriched over $\Zb$-graded complexes (over $k$) with differentials increasing degree. That is,
a category $\Cc$ is a dg-category if for $x$ and $y$ in $\Ob\Cc$
the set of morphsisms
\[\Cc(x,y)\]
forms a $\Zb$-graded complex of $k$-vector spaces. We will write
$(\Cc^\bullet(x,y),d)$ for this complex, if we need to reference the degree or differential in the complex. In addition, the
composition, for $x,y,z\in \Ob\Cc$
\[
\Cc(y,z)\otimes \Cc(x,y)\to
\Cc(x,z)\] is a morphism of complexes. Furthermore, there are obvious associativity and unit axioms. 
\end{defn}
\begin{defn}\label{Ho}
Given a dg-category $\Cc$, one can form the subcategory $Z^0\Cc$ which has the same objects as $\Cc$ and whose morphisms  from an object $x\in \Cc$ to an object $y\in \Cc$ are the degree $0$ closed morphisms in $\Cc(x,y)$.  The homotopy category $\Ho \Cc$ which has the same objects as $\Cc$ and whose morphisms  are the $0$th cohomology  
\[\Ho\Cc(x,y)=H^0(\Cc(x,y)).\]
A degree $0$ closed morphism between objects $x$ and $y$ in a dg-category $C$ is said to be a {\em homotopy equivalence} if it gives an isomorphism in the homotopy category.  
 \end{defn}
The dg-functors from a dg-category $C$ to the dg-category of chain complexes over $k$ themeselves form a dg-category, called  $\modl - \Cc$.  
The {\em Yondea embedding} of a dg-category $\Cc$ will refer to the functor 
\[\Cc \to \modl - \Cc
\]
given by 
\[c \mapsto \Cc(\cdot, c)
.\]
Let $\Cc$ and $\Dc$ be dg-categories.   A dg-morphism $F$  includes maps of chain complexes $\Cc(x,y) \to \Dc(F(x),F(y))$ for each $x,y \in \Cc$.  Notice that a dg-morphism $F$ takes homotopy equivalences to homotopy equivalences and also if $f$ is invertible in $\Ho\Cc$ then $F(f)$ is invertible in $\Ho\Dc$.  Following To\"{e}n \cite{T}, a dg-morphism $F:\Cc\to \Dc$ between two dg-categories is a {\em quasi-equivalence} if the induced morphisms $\Cc(x,y) \to \Dc(F(x),F(y))$ are quasi-isomorphisms for every pair of objects $x,y$ in $\Cc$ and also the induced functor $\Ho\Cc \to\Ho\Dc$ is essentially surjective. A dg-functor $F$ from $\Dc$ to $\modl-\Cc$ is called a {\it quasi-functor} from $\Dc$ to $\Cc$ if for each $d \in \Dc$, the object $F(d)$ is quasi-representable.  The following lemma consolidates the central results in the theory of dg-categories that we will use in the proof of our main theorem.  This is a consequence of a lemma of Keller (section 7.1 of \cite{Ke}).

\begin{lem}\label{Keller}
Suppose that   $\Cc$ and $\Dc$ are dg-categories, which are full dg-subcategories of dg-categories $C_{big}$ and $D_{big}$, so $\Cc \subset C_{big}$ and $\Dc \subset D_{big}$.    Let $F$ be a dg-functor from $C_{big}$ to $D_{big}$ which carries $C$ into $D$.   Let $G$ be a dg-functor from $D_{big}$ to $C_{big}$ which is right adjoint to $F$.  Say we have a collection of objects $c_d$ for each $d \in D$ and quasi-isomorphisms $h_{c_{d}} \to h_{G(d)}$  for each $d \in D$.  Let  $\Theta_d \in C(c_d,G(d))$ correspond to the identity of $c_d$ under $h_{c_{d}}(c_d) \to h_{G(d)}(c_d)$.  Assume also that that 
\[\Lambda_d= \eta_d \circ F(\Theta_d) :F(c_d) \to d
\] are homotopy equvialences in the category $D$, where $\eta:F \circ G \to \text{id} $ is the counit of the adjunction.  Suppose also that for each $c \in C,$ that $c$ and $c_{F(c)}$ are homotopy equivalent.  Then $F|_{C}$ is a dg-quasi-equivalence from $C$ to $D$. \end{lem}

{\bf Proof.}

Because dg-quasi-equivalence is an equivalence relation on dg-categories, it is enough to show that there exists a full dg sub-category $\iota: C_{small} \to C$ such that both of the arrows below are dg-quasi-equivalences 
\[C  \stackrel{\iota}\longleftarrow C_{small} \stackrel{F|_{ C_{small}}}\longrightarrow D
.\]
Let $C_{small}$ be the full dg subcategory with objects $\{c_d|d \in D\}$.  Then because each object $c$ in $C$ is homotopy equivalent to  $c_{F(c)}$, the inclusion $\iota$ of $C_{small}$ into $C$ is a dg-quasi-equivalence.  In order to show that $F|_{ C_{small}}$ is also a dg-quasi-equivalence using Keller's lemma in 7.2 of \cite{Ke}, we will need two ingredients: firstly the collection of objects 
\[\{(c_d, d)|d \in D\}  \subset \text{ob}(C_{small}) \times \text{ob}(D).
\] 
The projections of this collection to $C_{small}$ and to $D$ are 
clearly onto.   Secondly, we will need the $C_{small}-D-$bimodule defined 
by \[(c,d) \mapsto D(F(c),d).\]  Consider the maps
\begin{equation}\label{pleasework}C_{small}(c',c_d) \to D(F(c'),d)
\end{equation}
defined by 
\[g \mapsto \Lambda_d \circ F(g) = \eta_{d} \circ F(\Theta_d) \circ F(g) =  \eta_{d} \circ F(\Theta_d \circ g) 
\]
for arbitrary elements $c' \in C_{small}$.  By post-composing with $D_{big}(F(c'),d) \cong C_{big}(c',G(d))$ this map becomes simply 
\[C_{small}(c',c_d) \to C_{big}(c',G(d))
\]
\[g \mapsto \Theta_d \circ g
\]
which is clearly a quasi-isomorphism as it is just the given map $h_{c_{d}}(c') \to h_{G(d)}(c')$.  Hence, (\ref{pleasework}) is a quasi-isomorphism for any $c' \in C_{small}$.
The map 
\[
D(d,d') \to D(F(c_d),d')
\]
\[s \mapsto s \circ \Lambda_d
\]
is of course also a quasi-isomorphisms for any $d' \in D$.  Therefore, by Keller's Lemma in 7.2 of \cite{Ke} we see that the restriction of $F$ to $C_{small}$ is a dg-quasi-equivalence between $C_{small}$ and $D$.
\ \hfill $\Box$

\ \ 
\noindent
 A curved dga (see \cite{PP},  \cite{S}) is a triple
\[\As=(\A, d, c)\] where $\A= \bigoplus_{i=0}^{\infty}\A^{i}$ is a (non-negatively)
graded algebra over a field  $k$ of characteristic $0$, with a derivation
\[d:\A^\bullet\to \A^{\bullet+1}\]
which satisfies the usual graded Leibniz relation and
\[d^2(a)=[c,a]\]
where $c\in \A^2$ is a fixed element (the curvature) such that the Bianchi identity $dc=0$ holds.  Given any curved dga, we can consider the ``zero part"  \[\As^{0}=(\A^{0}, 0, 0)\] and conversely, any $k$-algebra $\A^{0}$ defines a curved dga in this way .

The category $\Pc_\As$ which was defined in \cite{B} consists of special types of $\As$-modules. Consider a $\Zb$-graded right module $\mE$ over $\A^0$.

 A
$\Zb$-connection $\Eb$  is a
$k$-linear  map
\[\Eb: \mE \otimes_{\A^0} \A \to \mE \otimes_{\A^0} \A
\]
of total degree one, which satisfies the usual Leibniz condition
\[  \Eb(e\omega)=(\Eb(e\otimes 1))\omega+(-1)^e e d\omega.\]

Such a connection is determined by its value on $\mE$. Let
$\Eb^k$ be the component of $\Eb$ which maps $\mE^\bullet$ to 
$\mE^{\bullet-k+1}\otimes_{\A^0}\A^k$.  It is clear that $\Eb^1$ is a
connection on each component $\mE^n$ in the ordinary sense when $n$ is even. When $n$ is odd it is a $\lambda$-connection for $\lambda=-1$. Also, notice that $\Eb^k$
is $\A^0$-linear for $k\ne 1$. 

\np Note that for a $\mathbb{Z}$-connection $\mathbb{E}$ on $\mE$ over a curved dga $\As=(\A, d, c)$,  the usual curvature $\mathbb{E}\circ\mathbb{E}$ is not $\A^0$-linear. However,  the {\em relative curvature}   
\[
F_\Eb(e)=\mathbb{E}\circ\mathbb{E}(e)+e\cdot c
\]
is $\A^0$-linear. 

For a curved dga $\As=(\A,d,c)$,  the dg-category $\Pc_{\As}$ was defined by:
\begin{enumerate}
\item An object $\tE=(\mE, \Eb)$ in
$\Pc_{\As}$, which we call a {\em cohesive module}, is a
$\Zb$-graded (but bounded in both directions) right module
$\mE = \bigoplus_{j=-\infty}^{\infty} \mE^j$ over $\A^0$ which is finitely generated and projective,
together with a $\Zb$-connection
\[\Eb: \mE \otimes_{\A^0}\A \to \mE \otimes_{\A^0}\A
\]
that satisfies the integrability condition that the relative curvature vanishes
\[F_\mathbb{E}(e)=\Eb\circ\Eb(e)+e\cdot
c=0\]
for all $e\in \mE$.

\item The morphisms of degree $k$,  ${\Pc^k_{\As}}(\tE_1,\tE_2)$ between two cohesive modules $\tE_1=(\mE_1,
\Eb_1)$ and $\tE_2=(\mE_2,\Eb_2)$ of degree $k$ are 
\[
\{\phi:\mE_1\otimes_{\A^0}
\A\to \mE_2\otimes_{\A^0} \A\,\,|\,\,\mbox{ of degree $k$ and }
\phi(ea)=\phi(e)a \,\,\, \forall a\in \A\}
\] with differential $d$ satisfying $d^2=0$ defined in the standard way
\[
d(\phi)(e)= \Eb_2(\phi(e))-(-1)^{\abs{\phi}}\phi(\Eb_1(e))
\]
Again, such a $\phi$ is determined by its restriction to
$\mE_1$ and if necessary we denote the component of $\phi$
that maps \begin{equation}\label{pcomp}\mE_1^\bullet\to
\mE_2^{\bullet+k-j}\otimes_{\A^0} \A^j \end{equation} by $\phi^j$. 
\end{enumerate}

Thus the morphisms in ${\Pc^k_{\As}}(\tE_1,\tE_2)$ are determined by maps in $\Hom^k_{\A^0}(\mE_1,\mE_2 \otimes_{\A^0} \A)$.
It is sometimes useful to consider the larger dg-category $\mathcal{C}_{\As}$: a {\em quasi-cohesive module } is the data of
$\tX=(\mX, \Xb)$ where $X$ is a
$\Zb$-graded  right module
over $\A$ 
together with a $\Zb$-connection
\[\Xb: X\otimes_{\A^0}\A \to X \otimes_{\A^0}\A
\]
that satisfies the integrability condition that the relative curvature 
\[F_{\Xb}(x)=\Xb \circ \Xb(x)+x\cdot
c=0\] for all $x\in \mX$. Thus, they differ from cohesive modules  by having possibly infinitely many nonzero
graded components as well as not being projective or finitely generated over $\A$.  Notice that $\Pc_{\As}$ is a full dg sub-category of $\mathcal{C}_{\As}$.

The shift functor on the category $\mathcal{C}_{\As}$ is defined as follows.  For 
$\tE=(\mE^\bullet, \mathbb{E})$ set 
$\tE[1]=(\mE[1]^\bullet,\mathbb{E}[1])$ where $\mE[1]^\bullet=\mE^{\bullet+1}$ and $\mathbb{E}[1]=-\mathbb{E}$.
It is easy to verify that $\tE[1]\in \mathcal{C}_{\As}$. Next for $\tE, \tF \in \mathcal{C}_{\As}$ and $\phi\in Z^0 \mathcal{C}_{\As}(\tE,\tF)$, define 
the cone of $\phi$, 
$\text{\bfseries\sf{Cone}}(\phi)=(Cone(\phi)^\bullet, \mathbb{C}_\phi)$ by
\[
Cone(\phi)^\bullet=\left(\begin{array}{c}\mF^\bullet \\  \oplus\\ \mE[1]^\bullet\end{array}\right)
\]
and 
\[
\mathbb{C}_\phi=\left(\begin{array}{cc} \mathbb{F} 		& \phi\\
                                 0              & \mathbb{E}[1] \end{array}\right)
\]
We then have a triangle of degree $0$ closed morphisms in  $\mathcal{C}_{\As}$
\begin{equation}\label{distinguishedtriangles}
\tE \stackrel{\phi}{\longrightarrow} \tF \longrightarrow \text{\bfseries\sf{Cone}}(\phi)    \longrightarrow \tE[1].  
\end{equation}
The cone satisfies the universal property of representing the functor
\[\mathcal{C}_{\As} \to C(k)
\]
given by 
\[\tM \mapsto Cone(\mathcal{C}_{\As}(\tM,\tE) \to \mathcal{C}_{\As}(\tM,\tF) )
.\]
When  $\tE, \tF \in \Pc_{\As}$ then clearly $\text{\bfseries\sf{Cone}}(\phi)\in \Pc_{\As}$.

Notice that for any curved dga $\Rs=(\mR,d_{\Rs}, c)$, there is a commutative diagram of dg-categories 
\begin{displaymath}
    \xymatrix{
       \Pc_{\Rs}   \ar[d] \ar[r] &   \Pc_{\Rs^{0}}  \ar[d]\\
             \mathcal{C}_{\Rs} \ar[r]                   &  \mathcal{C}_{\Rs^{0}} }
\end{displaymath}
where the horizonal arrows are defined on objects by 
\[\tE=(\mE,\mathbb{E}) \mapsto \tE^{0}= (\mE,\mathbb{E}^{0})
\] and on morphisms by 
\[\phi \mapsto \phi^{0}
.\]
 Let $\As$ be a curved dga.  Proposition 1 and Proposition 2 of \cite{BK} show that the dg-category $\Pc_{\As}$ is pretriangulated in the sense defined by Bondal and Kapranov, \cite{BK}. Therefore, the category $\Ho\Pc_{\As}$ is triangulated with the collection of distinguished triangles being isomorphic to those of the form \ref{distinguishedtriangles}.
The standard fully faithful Yoneda embedding is given by 
\[Z^0(\Pc_{\As}) \to \modl - \Pc_{\As}
\]

\[\tE \mapsto h_{\tE}= \Pc_{\As}(\cdot, \tE).
\]
while there is a similar functor 
\begin{equation}\label{quasiYoneda}Z^0(\mathcal{C}_{\As}) \to \modl - \Pc_{\As}
\end{equation}

\[\tX \mapsto \tilde{h}_{\tX}= \mathcal{C}_{\As}(\cdot, \tX).
\]
We need the following proposition from \cite{B}
\begin{prop}\label{characterization}
A closed morphism $\phi \in  \Pc^{0}_{\As}(\tE_1,\tE_2)$ is a homotopy equivalence if and only if $\phi^0:(\mE_1,\Eb_{1}^0) \to (\mE_2,\Eb_{2}^0)$ is a quasi-isomorphism of complexes of $\mA^0$ modules.
\end{prop}
We will also need the following theorem from \cite{B}
\begin{thm}\label{qf}
Suppose $\As=(\A, d, c)$ is a curved dga. Let $\tX=(\mX,\mathbb{X})$ be a quasi-cohesive module over $\As$. 
Then there is an object $\tH=(\mH, \mathbb{H})\in \Pc_\As$ such that  $\widetilde{h}_\tX$ is quasi-isomorphic to $h_{\tE}$; that is $\widetilde{h}_\tX$ is quasi-representable, under either of the two following conditions:
\begin{enumerate}
\item $X$ is a quasi-finite quasi-cohesive module.
\item $\A$ is flat over $\A^0$ and there is a bounded complex $(\mathcal{H},d_{\mH})$ of finitely generated projective right $\A^0$-modules and an $\A^0$-linear quasi-isomorphism $\Theta^0:(H,d_{H})\to (\mX,\mathbb{X}^0)$.
\end{enumerate}
In the second case we also have $\mathbb{H}^{0} = d_{H}$ and $H=\mH$.
\end{thm}
\begin{rem}\label{hope}The quasi-isomorphism $h_{\tH} \to \tilde{h}_{\tX}$ produced in the proof of (2) is of the form $\tilde{h}_{\Theta}$ for a specific element $\Theta \in \mathcal{C}_{\As}(\tH,\tX)$. Note also that one has in particular a quasi-isomorphism $\Pc_{\As}(\tH,\tH) \to \mathcal{C}_{\As}(\tE,\tX)$.  It can easily be seen that the zero component of the image of the identity $\text{id} \in \Pc_{\As}(\tH,\tH)$ inside $\mathcal{C}_{\As}(\tH,\tX)$ agrees with the original map of complexes $(H,d_{H})\to (\mX,\mathbb{X}^0)$. 
\end{rem}

\section{Pulling Back and Pushing Forward}
We now discuss a construction of functors between categories of the
form $\Pc_{\As}$. Given two curved dga's,
$\As_1=(\A_1,d_1,c_1)$ and $\As_2=(\A_2,d_2,c_2)$ a
homomorphism from $\As_1$ to $\As_2$ is a pair
$\tf=(f,\omega)$ where $f:\A_1 \to\A_2$ is a morphism
of graded algebras, $\omega\in \A_2^1$ and they satisfy
\begin{enumerate}
\item $f(d_1a_1)=d_2f(a_1)+[\omega,f(a_1)]$ and \item
$f(c_1)=c_2+d_2\omega+\omega^2$.
\end{enumerate}
Given a homomorphism of curved dga's $\tf=(f,\omega)$ we define a dg
functor
\[
\mathcal{C}_{\tA_{1}} \stackrel{\tf^{*}}\longrightarrow   \mathcal{C}_{\tA_{2}}\]
as follows. Given $\tE=(\mE,\Eb) \in \mathcal{C}_{\As_1}$ set $\tf^*(\tE)$ to be the cohesive module over
$\As_2$
\[(\mE \otimes_{\A^{0}_1} \A^{0}_2, \Eb_{\tA_2})\]
where $\Eb_{\tA_2}$ is determined by 
\[\Eb_{\tA_2}(e)=\Eb(e)+(-1)^{j}e\otimes \omega \] 
for $e \in \mE^{j}$.  One
checks that $\Eb_{\tA_2}$ is still a $\Eb$-connection and satisfies
\[ (\Eb_{\tA_2})^2(e\otimes b)=-(e\otimes b)c_2\]
The functor defined above takes objects in $\mathcal{P}_{\tA_{1}}$ to objects in $\mathcal{P}_{\tA_{2}}$.   Given $\psi \in \mathcal{C}_{\As_1}(\tE,\tF)$, the morphism $\tf^{*}\psi \in \mathcal{C}_{\As_1}(\tf^{*}\tE,\tf^{*}\tF)$  is determined by the composed map 
\[\mE \stackrel{\psi}\to \mF \otimes_{\mA^0_1} \mA_1 \to \mF \otimes_{\mA^0_1} \mA_2.
\]
Given composable morphisms $\tf, \tg$, there are natural equivalences $(\tf \circ \tg)^{*} \Longrightarrow \tg^{*} \circ \tf^{*}$ which satisfy the obvious coherence relation.   
\subsection{}

There is also a functor in the other direction 

\[ \modl - \Pc_{\As_{1}}
\stackrel{\tt f_{*}}\longleftarrow
 \modl - \Pc_{\As_{2}}.
\]
It is defined for any $\tM \in \mathcal{P}_{\tA_{1}}$ and $\tN \in \modl - \Pc_{\As_{2}}$  by 
\[(\tf_{*}\tN) (\tM) = \tN(\tf^{*}(\tM)).
\]

Suppose that $\tE$ and $\tF$ are in $\mathcal{P}_{\tA_{2}}$, and $\psi \in \modl - \Pc_{\As_{2}}(\tE,\tF)$.  We define $\tf_{*}\psi\in  \modl - \Pc_{\As_{1}}(\tf_{*}\tE,\tf_{*}\tF)$ by 

\[(\tf_{*}\psi)(\tM):(\tf_{*}\tE)(\tM) \to (\tf_{*}\tF)(\tM)
\]
\[g \to \psi(\tf^{*}\tM) \circ g.
\]
Given composable morphisms $\tf, \tg$, there are natural equivalences $(\tf \circ \tg)_{*} \Longrightarrow \tf_{*} \circ \tg_{*}$ which satisfy the obvious coherence relation.   
\subsection{}\label{FlatPush}Let $\tf=(f,\omega)$ be a map between $\tA_{1}=(\mA_{1},d_{1},c_{1})$ and  $\tA_{2}=(\mA_{2},d_{2},c_{2})$. 
Suppose now that we are in the special case that the natural map 
\[\mA^{0}_{2} \otimes_{\mA^{0}_{1}}  \mA_{1} \to \mA_{2}
\]
is an isomorphism.
Then we we will define a functor 
\begin{equation}\label{qPush}\mathcal{C}_{\tA_{1}}
\stackrel{\tt f_{*}}\longleftarrow
 \mathcal{C}_{\tA_{2}}.
\end{equation}
such that the diagram 
\[
\xymatrix{
\mathcal{C}_{\tA_{1}}  \ar[d]_{\tilde{h}} & \ar[l]_{\tf_{*}} \mathcal{C}_{\tA_{2}} \ar[d]^{\tilde{h}} \\
 \modl - \Pc_{\As_{1}} & \ar[l]_{\tf_{*}}  \modl - \Pc_{\As_{2}}}
\]
commutes where the vertical arrows are the functors defined in (\ref{quasiYoneda}).
Notice that in this case there is an isomorphism \begin{equation}\label{thiscase}\mM \otimes_{\mA^{0}_{1}} \mA_{1} \cong \mM \otimes_{\mA^{0}_{2}} \mA_{2}.\end{equation}
Therefore we can simply define
\[\tf_{*}(\mM,\Mb) = (\mM,\Mb_{\tf})
\]
where it understood that we use the isomorphism (\ref{thiscase}) in order to define the action of $\Mb$ on $\mM \otimes_{\mA^{0}_{1}} \mA_{1} $ and $\Mb_{\tf}$ is determined by 
\[\Mb_{\tf}(m) = \Mb(m)-(-1)^{j}m \otimes \omega,\] for $m \in \mathcal{M}^j$.
For any 
$\tM \in \mathcal{C}_{\tA_{1}}$ and $\tN \in \mathcal{C}_{\tA_{2}}$ we have natural adjointness isomorphisms 
\begin{equation}\label{adj}\mathcal{C}_{\tA_{2}}(\tf^{*}\tM,\tN) \cong \mathcal{C}_{\tA_{1}}(\tM,\tf_{*}\tN).
\end{equation}
Suppose that $\tE$ and $\tF$ are in $\mathcal{C}_{\tA_{2}}$, and $\beta \in \mathcal{C}_{\tA_{2}}(\tE,\tF)$.   We define $\tf_{*}\beta \in \mathcal{C}_{\tA_{1}}(\tf_{*}\tE,\tf_{*}\tF)$ to correspond via (\ref{adj}) to the element in  $\mathcal{C}_{\tA_{2}}(\tf^{*}\tf_{*}\tE,\tF)$ given by precomposing $\beta$ with the element of  $\mathcal{C}_{\tA_{2}}(\tf^{*}\tf_{*}\tE,\tE )$ which corresponds via (\ref{adj}) to the identity of $\mathcal{C}_{\tA_{2}}(\tf_{*}\tE,\tf_{*}\tE )$ .
\section{Homotopy Fiber Products}\label{HFP}
The purpose of this section is to explain the model we use for the homotopy fiber product of dg-categories and to describe three Mayer-Vietoris type consequences of a homotopy fiber product of dg-categories. These take the form long exact sequences of invariants associated to the four dg-categories appearing in the homotopy fiber product diagrams. The long exact sequences of invariants can in particular be applied to the homotopy fiber product diagram appearing in our main theorem \ref{main}.
Let B,C,D be dg-categories, along with dg-functors
\begin{equation}\label{setup}
\xymatrix{
 & B \ar[d]^{G} \\
C \ar[r]^{L} & D. }
\end{equation}
We need to use a suitable homotopy fiber product of the diagram (\ref{setup}) of dg-categories, where the notion of weak equivalence is given by dg-quasi-equivalence.  We chose the dg-category given by Drinfeld in section 15, Appendix IV (page 39) of \cite{Dr} and by Tabuada in Chapter 3, Definition 3.1 of \cite{Ta}.  Tabuada has constructed a model category structure on the category of dg-categories.  
Using this notion of homotopy fiber product of dg-categories, we consider the dg-category $B \times_{D}^{h} C$ with objects:
\begin{equation}
\begin{split}
\text{ob}(B \times_{D}^{h}C) =  \{ &(\tM \in B, \tN \in C,  \phi \in D^{0}(G(\tM) , L(\tN))) \\
& \text{such that} \phantom{x} \phi  \phantom{x} \text{is closed and becomes invertible in}  \phantom{x} H^{0}(D)\} 
\end{split}
\end{equation}

The morphisms are given by the complex which in degree $i$ is
\begin{equation}
\begin{split}
 & {(B \times_{D}^{h}C)}^{i}((\tM_{1},\tN_{1},\phi_{1}),(\tM_{2},\tN_{2},\phi_{2})) \\  &= B^i(\tM_{1} ,\tM_{2})\oplus  C^i( \tN_{1},\tN_{2}) \oplus D^{i-1}(G(\tM_{1}),L(\tN_{2}))
\end{split}
\end{equation}
with composition 
\begin{equation}\begin{split} & (B \times_{D}^{h}C)((\tM_2,\tN_2,\phi_2),(\tM_{3},\tN_{3},\phi_{3})) \otimes  (B \times_{D}^{h}C)((\tM_{1},\tN_{1},\phi_{1}),(\tM_2,\tN_2,\phi_2)) \\ & \to (B \times_{D}^{h}C)((\tM_{1},\tN_{1},\phi_{1}),(\tM_{3},\tN_{3},\phi_{3}))
\end{split}
\end{equation}
given by 
\[(\mu',\nu',\gamma')  (\mu,\nu,\tau) = (\mu' \mu, \nu' \nu, \gamma' G(\mu)+ L(\nu') \gamma ).
\]
The 
differential is given by 
\[d(\mu,\nu,\tau) = (d \mu,d \nu,d \gamma+ \phi_{2} G(\mu)- (-1)^{i}L(\nu) \phi_{1}).
\]
We will also need the dg sub category $B \stackrel{\to}\times_{D} C$ which has the same objects but 
\[(B \stackrel{\to}\times_{D} C)^{i}((\tM_{1},\tN_{1},\phi_{1}),(\tM_{2},\tN_{2},\phi_{2}))=\{(\mu,\nu) |\phi_{2} G(\mu)= (-1)^{i}L(\nu) \phi_{1}\}
\] 
where the differential and the composition are defined component-wise.  This dg-category itself has a full dg sub-category given by the actual fiber product $B \times_{D} C$ where $\phi$ is required to be an isomorphism.  In summary, we have dg-categories 
\[B \times_{D} C \subset B \stackrel{\to}\times_{D} C \subset B \times^{h}_{D} C
.\]

Notice that we can define  $B \times_{D}^{h}C$ bi-modules $\tilde{D}$ given by 
\[((\tM_{1},\tN_{1},\phi_{1}),(\tM_{2},\tN_{2},\phi_{2})) \mapsto D(G(\tM_{1}),L(\tN_2))[1]
\]
and $(B \times_{D}^{h}C,\tilde{D}) $ given by
\[((\tM_{1},\tN_{1},\phi_{1}),(\tM_{2},\tN_{2},\phi_{2})) \mapsto B(\tM_{1},\tM_2) \oplus C(\tN_{1},\tN_2).
\]
Together with $B \times_{D}^{h}C$ as a bi-module over itself, there is a short exact sequence of $B \times_{D}^{h}C$ bi-modules:
\[0\to \tilde{D} \to B \times_{D}^{h}C \to 
(B \times_{D}^{h}C,\tilde{D}) \to 0.\]  Therefore we have: 
\begin{lem}\label{les} Let B,C,D be dg-categories, along with dg-functors
as in equation (\ref{setup}). There is a long exact sequence of $H^{0}(B \times_{D}^{h}C)$ bi-modules

\[\xymatrix{
& &  \hspace{12pt}\cdots\hspace{12pt} \ar@{->}[r]  & H^{-1}(B \times_{D}^{h}C,\tilde{D})  \ar `r[d] `_l[lll] ^{} `^d[dlll] `^r[dll] [dll] \\
         & H^{0}(\tilde{D})  \ar@{->}[r]  &  H^{0}(B \times_{D}^{h}C)   \ar@{->}[r] &
                H^{0}(B \times_{D}^{h}C,\tilde{D})    \ar `r[d] `_l[lll] ^{} `^d[dlll] `^r[dll] [dll] \\
            &    H^{1}(\tilde{D})  \ar@{->}[r] &   \hspace{12pt}\cdots\hspace{12pt}   &
.}\]
There is also a short exact sequence of Hochschild cochain complexes of the dg-category $A=B \times_{D}^{h}C$ with coefficients in these three bi-modules.  It leads to (see \cite{LV}) a long exact sequence in Hochschild cohomology
\[0 \to HH^{0}(A,\tilde{D}) \to HH^{0}(A,A) \to HH^{0}(A,(B \times_{D}^{h}C,\tilde{D})) \to HH^{1}(A,\tilde{D}) \to \cdots
\]
\end{lem}
\begin{lem}\label{getback}Let B,C,D be dg-categories, along with dg-functors
as in equation (\ref{setup}). A closed morphism $(\mu,\nu,\tau)\in {(B \times_{D}^{h}C)}^{0}((\tM_{1},\tN_{1},\phi_{1}),(\tM_{2},\tN_{2},\phi_{2}))$ is a homotopy equivalence if and only if $\mu$ and $\nu$ are homotopy equivalences in  $B$ and $C$ respectively. \end{lem}

{\bf Proof.}  Suppose that  $(\mu,\nu,\tau)$ is a homotopy equivalence.
The necessity of $\mu$ and $\nu$ being homotopy equivalences follows from the fact that right composition with $(\mu,\nu,\tau)$ is a quasi-isomorphism,
\[ {(B \times_{D}^{h}C)}^{i}((\tM_2,\tN_2,\phi_2),(\tM_{3},\tN_{3},\phi_{3})) \to  {(B \times_{D}^{h}C)}^{i}((\tM_1,\tN_1,\phi_1),(\tM_{3},\tN_{3},\phi_{3}))
\]
given by 
\[(\mu',\nu',\gamma') \mapsto (\mu'\mu,\nu'\nu, \gamma G(\mu) + L(\nu') \gamma)
\]
 and hence, the composition with $\mu$ is a quasi-isomorphism 
 \[B(\tM_2,\tM_3) \to B(\tM_1,\tM_3)
 \]
\[\mu' \mapsto \mu' \mu
\]
 and similarly for $\nu$. Hence, $\mu$ and $\nu$ are homotopy equivalences.  Conversely, suppose that $\mu$ and $\nu$ are homotopy equivalences.
Consider a closed element in \\ ${(B \times_{D}^{h}C)}^{0}((\tM_{1},\tN_{1},\phi_{1}),(\tM_{2},\tN_{2},\phi_{2}))$  of the form
\[(\mu,\nu,\tau)\in   B^{0}(\tM_{1} ,\tM_{2})\oplus C^{0}( \tN_{1},\tN_{2}) \oplus D^{i-1}(G(\tM_{1}),L(\tN_{2})).
\]
Let us show that $(\mu,\nu,\tau)$ is invertible in the homotopy category.  Indeed for any $(\tM_{3},\tN_{3},\phi_{3}) \in B \times_{D}^{h}C$, right composition with $(\mu,\nu,\tau)$ induces a quasi-isomorphism 
\[\tilde{D}((\tM_2,\tN_2,\phi_2),(\tM_{3},\tN_{3},\phi_{3})) \to \tilde{D}((\tM_1,\tN_1,\phi_1),(\tM_{3},\tN_{3},\phi_{3}))
\] 
given by 
\[\gamma' \mapsto \gamma' G(\mu)
\]
and also a quasi-isomorphism 
\[(B \times_{D}^{h}C,\tilde{D})((\tM_2,\tN_2,\phi_2),(\tM_3,\tN_3,\phi_3)) \to (B \times_{D}^{h}C,\tilde{D})((\tM_1,\tN_1,\phi_1),(\tM_{3},\tN_{3},\phi_{3}))\]
given by 
\[(\mu',\nu') \mapsto (\mu' \mu, \nu' \nu)
.\]
Now an appeal to Lemma \ref{les} finishes the proof: the above maps induce isomorphisms on $H^{0}(\tilde{D}) $ and on $H^{0}(B \times_{D}^{h}C,\tilde{D})$ and hence on $H^{0}(B \times_{D}^{h}C)$.
\ \hfill $\Box$

\section{Patching of Complexes}

For a ring $\mathcal{R}$ let $P(\mathcal{R})$ be the category of projective modules over $\mathcal{R}$.   Let $P_{f}(\mathcal{R})$ be the category of finitely generated projective modules over $\mathcal{R}$. 
Let $\mA^{0},\mB^{0},\mC^{0},\mD^{0}$ be rings, along with a commutative diagram of ring homomorphisms
\begin{equation}\label{setup}
\xymatrix{
\mA^0 \ar[r]^{f^0} \ar[d]^{k^0} & \mB^0 \ar[d]^{g^0} \\
\mC^0 \ar[r]^{l^0} & \mD^0 }.
\end{equation}
such that
\begin{ass}\label{MilnorConditions}
\ \ 
\begin{enumerate}
\item the maps $f^{0},g^{0},k^{0},l^{0}$ form a fiber product:
 \[\mA^0 \cong \mB^0 \times_{\mD^0} \mC^0 \]
\item $l^{0}:\mathcal{C}^{0} \to \mathcal{D}^{0}$ is surjective. 
 \end{enumerate}
\end{ass}
\subsection{Categories of Projective Modules}\label{CatProjMod}
Recall that the fiber product category $P(\mB^{0}) \times_{P(\mD^{0})} P(\mC^{0})$ has objects consisting of a pair of objects $\mM \in P(\mB^{0})$ and $\mN \in P(\mC^0)$ together with an isomorphism $\alpha: \mM \otimes_{\mB^0} \mD^0 \to \mN \otimes_{\mC^0} \mD^0$ while a morphism from $(\mM_1,\mN_1,\alpha_1)$ to $(\mM_2,\mN_2,\alpha_2)$ consists of a pair of morphisms $u \in  P(\mB^{0})(\mM_1,\mM_2)$ and $v \in  P(\mC^{0})(\mN_1,\mN_2)$ whose images in $P(\mD^0)(\mM_1 \otimes_{\mB^0} \mD^0, \mM_2 \otimes_{\mB^0} \mD^0)$ and $P(\mD^0)(\mN_1 \otimes_{\mC^0} \mD^0, \mN_2 \otimes_{\mC^0} \mD^0)$  intertwine $\alpha_{1}$ and $\alpha_{2}$.
Under the assumptions, \ref{MilnorConditions} Milnor's theorem from the second chapter of \cite{M} establishes an equivalence of categories 
\[R: P(\mA^{0})  \longrightarrow P(\mB^{0}) \times_{P(\mD^{0})} P(\mC^{0}) 
\]
\[R(\mE) = (\mE \otimes_{\mA^{0}}\mB^{0}, \mE \otimes_{\mA^{0}}\mC^{0}, \kappa)
\]
with inverse 
\begin{equation}\label{MilnorPsi} P(\mA^{0})  \longleftarrow P(\mB^{0}) \times_{P(\mD^{0})} P(\mC^{0})  :\psi
\end{equation}

\[\psi(\mM,\mN,\phi) = \ker[\mM \oplus \mN \to \mN\otimes_{\mC^{0}} \mD^{0} ].
\]
In fact, the functor $R$ is a left adjoint to the functor $\psi$ and the unit and counit adjuctions are isomorphisms for every object.  The same functors respect the property of being finitely generated and so establish the equivalence of categories \[P_{f}(\mA^{0})  \longrightarrow P_{f}(\mB^{0}) \times_{P_{f}(\mD^{0})} P_{f}(\mC^{0}) 
.\]
\subsection{Categories of Complexes of Modules}

Consider the category $\mathcal{C}_{\tB^0} \stackrel{\to}\times_{\mathcal{C}_{\tD^0}} \mathcal{C}_{\tC^0}$.  The objects are chain complexes $(\mathcal{M},d_{\mM})$ over $\mB^{0}$, $(\mathcal{N},d_{\mN})$ over $\mC^{0}$ and a degree zero morphism \[\alpha: \mathcal{M} \otimes_{\mB^{0}}\mD^{0} \to \mathcal{N} \otimes_{\mC^{0}}\mD^{0}
\] of complexes of $\mD^{0}$ modules commuting ($\alpha \circ (d_{\mM} \otimes \text{id}) = (d_{\mN} \otimes \text{id}) \circ \alpha$) with the differentials.   There is a dg-functor given by restriction 
\[R:\mC_{\tA^0} \to \mathcal{C}_{\tB^0} \stackrel{\to}\times_{\mathcal{C}_{\tD^0}} \mathcal{C}_{\tC^0}
\]
which sends 
\[(\mathcal{E},d_{\mathcal{E}}) \to (\tf^{0*}(\mathcal{E},d_{\mathcal{E}}), \tk^{0*}(\mathcal{E},d_{\mathcal{E}}),\kappa)
.\]
For any object $((\mathcal{M},d_{\mM}), (\mathcal{N},d_{\mN}),\alpha)\in \mathcal{C}_{\tB^0} \stackrel{\to}\times_{\mathcal{C}_{\tD^0}} \mathcal{C}_{\tC^0}$, consider the difference of the maps 
\[\mathcal{M} \to \mathcal{M} \otimes_{\mB^{0}}\mD^{0} \stackrel{\alpha}\to \mathcal{N} \otimes_{\mC^{0}}\mD^{0}
\]
and 
\[\mN \to \mN \otimes_{\mC^{0}}\mD^{0}
\]
as a single map 
\[
\mathcal{M} \oplus \mathcal{N} \to  \mathcal{N} \otimes_{\mC^{0}}\mD^{0}.
\]
It clearly intertwines the differential $d_{\mM} \oplus d_{\mN}$ with the differential $d_{\mathcal{N}} \otimes \text{id}_{\mD^0}$.  We can use 
this to define a dg-functor 
\begin{equation}\label{psidef}  \mC_{\tA^0} \longleftarrow \mathcal{C}_{\tB^0} \stackrel{\to}\times_{\mathcal{C}_{\tD^0}} \mathcal{C}_{\tC^0}: \psi
\end{equation}
which assigns 
\[\psi(\mM,d_{\mM},\mN,d_{\mN},\alpha) =(\ker[\mathcal{M} \oplus \mathcal{N} \to  \mathcal{N} \otimes_{\mC^{0}}\mD^{0}], d_{\mM} \oplus d_{\mN}).
\]
Notice that 
\[\psi(\mM,d_{\mM},\mN,d_{\mN},\alpha)^{i} = \ker[\mathcal{M}^{i} \oplus \mathcal{N}^{i} \to  \mathcal{N}^{i} \otimes_{\mC^{0}}\mD^{0}] = \psi(\mM^i,\mN^i),
\] 
so in each component, this is just the original map \ref{MilnorPsi}.  
\begin{lem}\label{Landsburg}
It was shown by Landsburg in \cite{L} that assuming conditions (1) and (2) from \ref{MilnorConditions} that $\psi$ sends the full subcategory of gluing data for complexes of projective modules where the morphism $\alpha$ is an isomorphism into the category of complexes of projective modules.  In particular we have a commutative diagram 
\begin{equation}\label{LandsburgDiagram2}
\begin{CD} \mathcal{P}_{\tA^0} 	@<<\psi<    \mathcal{P}_{\tB^0}  \times_{ \mathcal{P}_{\tD^0} }   \mathcal{P}_{\tC^0}  \\ @VVV @VVV\\ 	\mathcal{C}_{\tA^0} @<<\psi< \mathcal{P}_{\tB^0}  \stackrel{\to}\times_{ \mathcal{P}_{\tD^0} }   \mathcal{P}_{\tC^0}. \end{CD}
\end{equation}
Therefore, the restriction map $R$ gives an equivalence of dg-categories 
\[\mathcal{P}_{\tA^0} \cong \mathcal{P}_{\tB^0} \times_{ \mathcal{P}_{\tD^0}} \mathcal{P}_{\tC^0}
\] with inverse $\psi$.
\end{lem}

\subsection{Enhanced categories of complexes of modules}

Consider the dg-category given by $\mathcal{C}_{\tB^0} \stackrel{h}\times_{\mathcal{C}_{\tD^0}} \mathcal{C}_{\tC^0}$.

There is the functor 
\begin{equation}\label{psitildedef}  \mathcal{C}_{\tA^0} \longleftarrow \mathcal{C}_{\tB^0} \times_{\mathcal{C}_{\tD^0}}^{h} \mathcal{C}_{\tC^0} :\tilde{\psi}
\end{equation} 
which assigns 
\[\tilde{\psi}(\mM,d_{\mM},\mN,d_{\mN},\alpha) =\text{Cone}[\mathcal{M} \oplus \mathcal{N} \to  \mathcal{N} \otimes_{\mC^{0}}\mD^{0}][-1],
\]
the cone is taken in the dg-category $\mathcal{C}_{\tA^0}$.
The functor $\tilde{\psi}$ is a right adjoint to 
\[R: \mathcal{C}_{\tA^0} \longrightarrow \mathcal{C}_{\tB^0} \times_{\mathcal{C}_{\tD^0}}^{h} \mathcal{C}_{\tC^0}.
\]
Landsburg observed in 2.4 of \cite{L} that $\psi$ and $\tilde{\psi}$ are quasi-isomorphic functors  when restricted to $\mathcal{P}_{\tB^0}  \stackrel{\to}\times_{ \mathcal{P}_{\tD^0} }   \mathcal{P}_{\tC^0} $.   In particular, we will usually consider the restriction of $\tilde{\psi}$ to the the subcategory 
\[ \mathcal{P}_{\tB^0} \times^{h}_{\mathcal{P}_{\tD^0}}  \mathcal{P}_{\tC^0}
 \subset \mathcal{C}_{\tB^0} \times_{\mathcal{C}_{\tD^0}}^{h} \mathcal{C}_{\tC^0}
\]
What is important for us is that the diagram

\begin{equation}\label{LandsburgDiagram}
\begin{CD} \mathcal{C}_{\tA^0} 	@<<\psi<    \mathcal{P}_{\tB^0}  \stackrel{\to}\times_{ \mathcal{P}_{\tD^0} }   \mathcal{P}_{\tC^0}  \\ @V=VV @VVV\\ 	\mathcal{C}_{\tA^0} @<<\tilde{\psi}< \mathcal{P}_{\tB^0}  \times^{h}_{ \mathcal{P}_{\tD^0} }   \mathcal{P}_{\tC^0} \end{CD}
\end{equation}
commute up to a natural transformation 
\begin{equation}\label{T}T: \psi \to \tilde{\psi}
\end{equation} which for each object 
of $\mathcal{P}_{\tB^0}  \stackrel{\to}\times_{ \mathcal{P}_{\tD^0} }   \mathcal{P}_{\tC^0}$, gives a quasi-isomorphism between the object's two images in $\mathcal{C}_{\tA^0}$.  This natural transformation is described more explicitly in the proof of Lemma \ref{DefCollection}.
The diagrams

\begin{equation}\label{zeroDiagOne}
    \xymatrix{
        R\circ \psi  \ar[dr] \ar[r]^{R\circ T} & R \circ \tilde{\psi}  \ar[d]\\
                               & \iota }
\end{equation}
and 
\begin{equation}\label{zeroDiagTwo}
    \xymatrix{
        \iota  \ar[dr] \ar[r] & \psi \circ R  \ar[d]^{T \circ R}\\
                               & \tilde{\psi} \circ R }
\end{equation}
of functors $ \mathcal{P}_{\tB^0} \stackrel{\to}\times_{\mathcal{P}_{\tD^0}}  \mathcal{P}_{\tC^0} \to \mathcal{C}_{\tB^0} \times_{\mathcal{C}_{\tD^0}}^{h} \mathcal{C}_{\tC^0}$ and   $ \mathcal{P}_{\tA^0} \to \mathcal{C}_{\tA^0}$ respectively commute up to quasi-isomorphism where $\iota$ is the natural inclusion.  The vertical arrow in the first diagram and the diagonal arrow in the second diagram are actually natural equivalences.  Since the horizontal arrow in the first diagram and the vertical arrow in the second diagram give a quasi-isomorphism for each object, we can conclude the same about the vertical arrow in the first diagram and the diagonal arrow in the second diagram.


\section{Descent}

Let ${\tA}=(\mathcal{A},d_{\tA},c_{\tA}),{\tB}=(\mathcal{B},d_{\tB},c_{\tB}),{\tC}=(\mathcal{C},d_{\tC},c_{\tC}),{\tD}=(\mathcal{D},d_{\tD},c_{\tD})$ be curved dg-algebras, along with degree $0$ maps of curved dg-algebras forming a commutative diagram:
\begin{equation}\label{setup2}
\xymatrix{
\tA \ar[r]^{\tf} \ar[d]^{\tk} & \tB \ar[d]^{\tg} \\
\tC \ar[r]^{\tl} & \tD. }
\end{equation}
and satisfying assumptions (\ref{MilnorConditions}).
We will primarily be interested in the homotopy fiber product of dg-categories $\mathcal{P}_{\tB} \times_{\mathcal{P}_{\tD}}^{h} \mathcal{P}_{\tC}$ corresponding to the diagram

\begin{equation}
\xymatrix{
 & \mathcal{P}_{\tB} \ar[d]^{\tg^{*}} \\
\mathcal{P}_{\tC} \ar[r]^{\tl^{*}} & \mathcal{P}_{\tD}. }
\end{equation}
There is a dg-functor:
\begin{equation}\label{RFirst}\mathcal{P}_{\tA} \stackrel{R}\longrightarrow \mathcal{P}_{\tB} \times_{\mathcal{P}_{\tD}}^{h} \mathcal{P}_{\tC}.
\end{equation}
\[R(\tS)= (\tf^{*}\tS, \tk^{*}\tS, \kappa)
\]
where $\kappa$ is the canonical isomorphism $\tg^{*}\tf^{*}\tS \to \tl^{*}\tk^{*}\tS $.
In order to define a reasonable functor in the opposite direction (but to $\mathcal{C}_{\tA}$) we will need an extra assumption.
\begin{ass}\label{flatpush}
 Assume that the natural maps induced by $\tf, \tg, \tk, \tl$ induce isomorphisms $\mB^{0} \otimes_{\mA^{0}} \mA\to \mB$,  $\mC^{0} \otimes_{\mA^{0}} \mA \to \mC$,  $\mD^{0} \otimes_{\mB^{0}} \mB \to \mD$,  $\mD^{0} \otimes_{\mC^{0}} \mC  \to \mD$.  
 \end{ass}
 Lets start with $(\tM, \tN, \phi) \in \mathcal{P}_{\tB} \times_{\mathcal{P}_{\tD}}^{h} \mathcal{P}_{\tC}$.  Due to our assumptions on $\tf, \tg, \tk, \tl$ we can define pushforward functors as were described in \ref{qPush}.
Consider the map $\lambda \in \mathcal{C}_\tA(\tf_{*}\tM\oplus\tk_{*}\tN ,\tf_{*} \tg_{*}\tl^{*} \tN)$ 
\begin{equation}\label{LambdaEq}\tf_{*}\tM \oplus \tk_{*}\tN \stackrel{\lambda} \longrightarrow \tf_{*}\tg_{*}\tl^{*}\tN 
\end{equation}
defined as the difference of the maps 
\[\tf_{*}\tM \longrightarrow \tf_{*}\tg_{*}\tg^{*}\tM \stackrel{ \tf_{*}\tg_{*} \phi}\longrightarrow \tf_{*}\tg_{*}\tl^{*}\tN 
\]
and 
\[\tk_{*}\tN \longrightarrow \tk_{*}\tl_{*}\tl^{*}\tN =  \tf_{*}\tg_{*}\tl^{*}\tN.
\]

Define $\tt{X}=(\mX,\Xb) \in \mathcal{C}_\tA$ as the cone of $\lambda$.  
\begin{defn}\label{ATilde}
Consider curved dg-algebras and degree $0$ maps of curved dg-algebras as in (\ref{setup2}) satisfying assumption \ref{MilnorConditions} and assumption \ref{flatpush}.  We now define a functor 
\[
 \mathcal{C}_{\tA} \stackrel{\tilde{A}}\longleftarrow \mathcal{C}_{\tB} \times_{\mathcal{C}_{\tD}}^{h} \mathcal{C}_{\tC}.
\]
In terms of objects, we send $(\tM, \tN, \phi)$ to 
\[\tilde{A}((\tM, \tN, \phi))= \tX= \text{\bfseries\sf{Cone}}(\lambda)\] 
where $\lambda$ was defined in (\ref{LambdaEq}).
 Suppose that $(\mu,\nu,\tau)$ is a morphism in $\mathcal{C}_{\tB} \times_{\mathcal{C}_{\tD}}^{h} \mathcal{C}_{\tC}$ from $(\tM_1, \tN_1, \phi_1)$ to $(\tM_2, \tN_2, \phi_2)$.  We need to define a map in $\mC_{\tA}(\tX_{1},\tX_{2})$ in order to define $\tilde{A}$ on morphisms.  We use the map
\[\left(\begin{array}{cc} \tf_{*}\tg_{*}\tl^{*}\nu		& \epsilon \\
                                 0              &  \tf_{*}(\mu) \oplus \tk_{*}(\nu) \end{array}\right): \mX_{1} \otimes_{\mA^0} \mA \to \mX_{2} \otimes_{\mA^0} \mA
\]
where $\epsilon$ is the map 
\[\tf_{*}\tM_{1} \oplus \tk_{*}\tN_{1} \to \tf_{*}\tM_{1} \to \tf_{*}\tg_{*}\tg^{*}\tM_{1} \stackrel{ \tf_{*}\tg_{*}\gamma}\to  \tf_{*}\tg_{*}\tl^{*}\tN_{2}\]
and we have used the presentations 
\[\mX_{i}= \left(\begin{array}{c}\mN_{i} \otimes_{\mC^0} \mD^0[-1] \\  \oplus\\ \mM_{i} \oplus \mN_{i}.\end{array}\right).
\]
\end{defn}
Notice that 
\begin{equation}\label{ADJ}\begin{split} 
 {(\mathcal{C}_{\tB} \times_{\mathcal{C}_{\tD}}^{h} \mathcal{C}_{\tC})}^{i}(R(\tS),(\tM,\tN,\phi))
&=  {(\mathcal{C}_{\tB} \times_{\mathcal{C}_{\tD}}^{h} \mathcal{C}_{\tC})}^{i}((\tf^{*}\tS, \tk^{*}\tS, c),(\tM,\tN,\phi)) \\ & =  \mathcal{C}^{i}_{\tB}(\tf^{*}\tS ,\tM)\oplus \mathcal{C}^{i}_{\tC}(\tk^{*}\tS,\tN) \oplus \mathcal{C}^{i-1}_{\tD}(\tg^{*}\tf^{*}\tS,\tl^{*}\tN) \\
&  \cong \mathcal{C}_{\tA}^{i}(\tS,\tf_{*}\tM \oplus \tk_{*}\tN \oplus \tf_{*}\tg_{*}\tl^{*}\tN[1]) \\
&=  \mathcal{C}_{\tA}^{i}(\tS,\tilde{A}(\tM,\tN,\phi))
.\end{split}
\end{equation}
This correspondence is natural in the variable $\tS$ and also in the variable $(\tM,\tN,\phi)$.  It intertwines the differentials.  Therefore, $\tilde{A}$ is a right adjoint to $R$.

Notice that if $\tX=\tilde{A}(\tM,\tN,\phi))$ the complex $(\mX,\Xb^{0})$ is none other than $\tilde{\psi}(\tM^0,\tN^0)$, where $\tilde{\psi}$ is the functor defined in (\ref{psitildedef}).  Using the natural adjunction $\eta_{\tE}: \tE \to \tilde{A}(R(\tE))$ we find 
\begin{lem}\label{natZero} Consider curved dg-algebras and degree $0$ maps of curved dg-algebras as in (\ref{setup2}) satisfying assumption \ref{MilnorConditions} and assumption \ref{flatpush}.   Let $\tE \in  \mathcal{P}_{\tA}$ be any object.  There is a degree zero, closed morphism in $\mC_{\tA}$  \[\eta_{\tE}: \tE \to \tilde{A}(R(\tE))\] such that the map \[\eta_{\tE}^{0}: \tE^0 \to \tilde{A}(R(\tE))^0 = \tilde{\psi}(R(\tE^0))\] (which is a morphism in $\mC_{\tA^0}$) is a quasi-isomorphism of complexes of $\mA^{0}$ modules.
\end{lem}
{\bf Proof.}
The map $\eta_{\tE}^{0}$ is simply the diagonal arrow in the diagram (\ref{zeroDiagTwo}) applied to $\tE^0$.  As discussed below that diagram , this is a quasi-isomorphism.
\ \hfill $\Box$

\begin{lem}\label{DefCollection}Consider curved dg-algebras and degree $0$ maps of curved dg-algebras as in (\ref{setup2}) satisfying assumption \ref{MilnorConditions} and assumption \ref{flatpush}.   Assume furthur that $\A$ is flat over $\A^0$.  Then for every $(\tM,\tN,\phi) \in \mathcal{P}_{\tB} \times_{\mathcal{P}_{\tD}}^{h} \mathcal{P}_{\tC}$, there is an object $\tH_{(\tM,\tN,\phi)} \in \mathcal{P}_{\tA}$ and a morphism $\Theta_{(\tM,\tN,\phi)} \in \mathcal{C}_{\tA}(\tH_{(\tM,\tN,\phi)},\tilde{A}(\tM,\tN,\phi))$ such that $\tilde{h}_{\Theta_{(\tM,\tN,\phi)}}$  is quasi-isomorphism of $\mathcal{P}_{\tA}$ modules
\[\tilde{h}_{\Theta_{(\tM,\tN,\phi)}}: h_{\tH_{(\tM,\tN,\phi)}} \to \tilde{h}_{\tilde{A}(\tM,\tN,\phi)}
\]
\end{lem}

{\bf Proof.}
Consider $\tX= \tilde{A}(\tM,\tN,\phi) \in \mathcal{C}_{\tA}$. 
Notice that from the definition of $\tX=(\mX,\Xb)$, it follows that 
\[\mX= \left(\begin{array}{c}\mN \otimes_{\mC^0} \mD^0[-1] \\  \oplus\\ \mM \oplus \mN \end{array}\right)
\] and that 
\[\mathbb{X}^{0}= \left(\begin{array}{cc} \mathbb{N}^{0}_{\tt D}[-1] 		& \lambda^{0} \\
                                 0              & \Mb^0 \oplus \Nb^0 \end{array}\right).
\]

Because the map $\mC^0 \to \mD^0$ is surjective, we get that $\mN \to \mN \otimes_{\mC^0} \mD^0$ is surjective, therefore the map $\alpha^0$ is as well and so the cohomology of $(\mX,\Xb^0)$ is generated by the second row, in other words by elements of $\mM \oplus \mN$.   If we consider the zeroth component of the equation $\phi \circ \mathbb{M}_{\tD} = \mathbb{N}_{\tD} \circ \phi$ we learn that 
\[\phi^{0} \circ \mathbb{M}^{0}_{\tD} = \mathbb{N}^{0}_{\tD} \circ \phi^{0}
.\]
Therefore, we can consider the complex of $\mA^{0}$ modules
\begin{equation}\label{Gdef}(\mG,d_{\mG}) = \psi(\tM^0,\tN^0,\phi^{0}) \in \mathcal{C}_{\tA^0}.
\end{equation}
The map 
\[T(\tM^0,\tN^0,\phi^{0}): (\mG,d_{\mathcal{G}}) \to (\mX,{\Xb}^0)=\tX^0
\] which was mentioned in (\ref{T}) simply includes the complex $\mG \subset \mM \oplus \mN$ defined in equation (\ref{Gdef}) into the second row of $\mX$ and is a quasiisomorphism of complexes of $\mA^0$ modules. 
Since $\phi^{0}$ need not be an isomorphism, we cannot apply Milnor's theorem from \cite{M} and conclude that $\mG$ is a projective $\mA^{0}$ module.  However, Landsburg actually shows in Lemma 1.4, page 362 and in the theorem on page 269 of \cite{L} that there exist bounded acyclic complexes of finitely generated and projective $\mB^{0}$ and $\mC^{0}$ modules $(\tilde{\mM},d_{\tilde{\mM}})$ and $(\tilde{\mN},d_{\tilde{\mN}})$, an {\it isomorphism} of complexes
\[\alpha: (\mM \oplus \tilde{\mM})\otimes_{\mB^0} \mD^0 \to (\mN\oplus  \tilde{\mN})\otimes_{\mC^0} \mD^0,
\] and a quasi-isomorphism 
\begin{equation}\label{firstquasi}(\mH,d_{\mH})=\psi(\mM \oplus \tilde{\mM},\Mb^{0} +d_{ \tilde{\mM}} ,\mN \oplus \tilde{\mN},\Nb^{0}+d_{ \tilde{\mN}},\alpha) \to (\mG,d_{\mG}).
\end{equation}
Note that Landsburg's Lemma which we gave in Lemma \ref{Landsburg} tells us that $(\mH,d_{\mathcal{H}})$ is a bounded complex of finitely generated projective $\mathcal{A}^{0}$ modules and so $(\mH,d_{\mathcal{H}})\in \mathcal{P}_{\tA^0}$.  
By combining this with (\ref{firstquasi}) we have a sequence of quasiisomorphisms of complexes of $\mA^{0}$ modules
\[(\mathcal{H},d_{\mathcal{H}}) \to (\mathcal{G},d_{\mathcal{G}}) \to (\mX,\Xb^{0})=\tX^0.
\]
Therefore by Theorem \ref{qf}  there exists a superconnection $\mathbb{H}$ on $\mH$ such that $\mathbb{H}^{0}=d_{\mathcal{H}}$ and ${\tt H}=(\mathcal{H},\mathbb{H}) \in \mathcal{P}_\tA$, in addition there is a map $\Theta_{(\tM,\tN,\phi)} \in \mathcal{C}_{\tA}(\tH,\tX)$ such that $\tilde{h}_{\Theta_{(\tM,\tN,\phi)}}:h_{\tH} \to \tilde{h}_\tX$ is a quasi-isomorphism.  In other words, $\tilde{h}_\tX$ is quasi-representable.

\ \hfill $\Box$
\begin{lem} \label{injectivity} Consider curved dg-algebras and degree $0$ maps of curved dg-algebras as in (\ref{setup2}) satisfying assumption \ref{MilnorConditions} and assumption \ref{flatpush}.  For every $\tE \in \mathcal{P}_{\tA}$,  $\tE$ and  $H_{R(\tE)}$ are homotopy equivalent.
\end{lem}
{\bf Proof.}
Recall from Lemmas \ref{natZero} and \ref{DefCollection} we have morphisms in $\mathcal{C}_{\tA}$, of the form 
\[\tE \stackrel{\eta_{\tE}}\longrightarrow \tilde{A}(R(\tE)) \stackrel{\Theta_{R(\tE)}} \longleftarrow \tH_{R(\tE)}
\]
such that the right arrow induces a quasi-isomorphism $\tilde{h}_{\tilde{A}(R(\tE))} \longleftarrow h_{\tH_{R(\tE)}}$.
Choose a closed, degree zero element  $t_{\tE} \in \mathcal{C}_{\tA}(\tilde{A}(R(\tE)) ,\tH_{R(E)})$ which maps to the identity under the map 
\[H^{0}(\tilde{A}(R(\tE)),\tH_{R(\tE)}) \to H^{0}(\tH_{R(\tE)},\tH_{R(\tE)})\]
given by $(\;) \circ \Theta_{R(\tE)}$.  Now $t_{\tE} \circ \Theta_{R(\tE)}$ is the identity in the homotopy category and therefore, since $\Theta_{R(\tE)}^0$ is a quasi-isomorphism, $t_{\tE}^0$ must be as well.   We also know that $\eta_{\tE}^{0}$ is a quasi-isomorphism from Lemma \ref{natZero}.  Consider 
\[\tE \stackrel{\eta_{\tE}}\longrightarrow \tilde{A}(R(\tE)) \stackrel{t_{\tE}}\longrightarrow \tH_{R(\tE)}.
\]
We know that $(t_{\tE} \circ \eta_{\tE})^{0}$ is a quasi-isomorphism and hence by Proposition \ref{characterization},  $t_{\tE} \circ \eta_{\tE} \in \mathcal{P}_{\tA}(\tE,  \tH_{R(\tE)})$ is a homotopy equivalence.
\ \hfill $\Box$

We now state the assumptions on the commutative diagram (\ref{setup}) of curved dgas and morphisms between them which will be used in our main theorem.
\begin{ass}\label{conditions}
\ \ 
\begin{enumerate}
\item the maps $f^{0},g^{0},k^{0},l^{0}$ form a fiber product in degree zero:
 \[\mA^0 \cong \mB^0 \times_{\mD^0} \mC^0 \]
\item $l^{0}:\mathcal{C}^{0} \to \mathcal{D}^{0}$ is surjective 
\item
$\A$ is flat over $\A^0$ 

 \item
  the natural maps induced by $\tf, \tg, \tk, \tl$ induce isomorphisms $\mB^{0} \otimes_{\mA^{0}} \mA\to \mB$,  $\mC^{0} \otimes_{\mA^{0}} \mA \to \mC$,  $\mD^{0} \otimes_{\mB^{0}} \mB \to \mD$,  $\mD^{0} \otimes_{\mC^{0}} \mC  \to \mD$.
 \end{enumerate}
\end{ass}
\begin{thm}\label{main}  Let $\tA,\tB,\tC,\tD$ be curved dg-algebra  with morphisms between them as in the diagram (\ref{setup}).  Suppose that these morphisms satisfy Assumptions \ref{conditions}.  Then the restriction functor $R$ is a dg-quasi-equivalence of dg-categories between $\mathcal{P}_{\tA}$ and $\mathcal{P}_{\tB} \times_{\mathcal{P}_{\tD}}^{h} \mathcal{P}_{\tC}$.
\end{thm}

{\bf Proof.}
We wish to apply Lemma \ref{Keller} to the dg-functor given by restriction 
\begin{equation}\mathcal{C}_{\tA} \stackrel{R}\longrightarrow \mathcal{C}_{\tB} \times_{\mathcal{C}_{\tD}}^{h} \mathcal{C}_{\tC}, 
\end{equation}
\[R(\tS)= (\tf^{*}\tS, \tk^{*}\tS, \kappa)
\]
where $\kappa$ is the canonical isomorphism $\tg^{*}\tf^{*}\tS \to \tl^{*}\tk^{*}\tS$, 
and its right adjoint (see Definition \ref{ATilde})
\begin{equation}\mathcal{C}_{\tA} \stackrel{\tilde{A}}\longleftarrow \mathcal{C}_{\tB} \times_{\mathcal{C}_{\tD}}^{h} \mathcal{C}_{\tC}, 
\end{equation}
in order to show that the restriction of $R$:
\begin{equation}\label{beEquiv}\mathcal{P}_{\tA} \stackrel{R}\longrightarrow \mathcal{P}_{\tB} \times_{\mathcal{P}_{\tD}}^{h} \mathcal{P}_{\tC}, 
\end{equation}
is a dg-quasi-equivalence.
  For the collection of objects required by Lemma \ref{Keller} we will use the collection 
\[\{ \tH_{(\tM,\tN,\phi)} \in \mathcal{P}_{\tA} |(\tM,\tN,\phi) \in  \mathcal{P}_{\tB} \times_{\mathcal{P}_{\tD}}^{h} \mathcal{P}_{\tC} \} 
\]
which were defined in Lemma \ref{DefCollection}.  
Consider an arbitrary object $(\tM,\tN,\phi) \in \mathcal{P}_{\tB} \times_{\mathcal{P}_{\tD}}^{h} \mathcal{P}_{\tC}$.
There is an element in 
\[\Theta_{(\tM,\tN,\phi)} \in \mathcal{C}_\tA(\tH_{(\tM,\tN,\phi)},\tilde{A}(\tM,\tN,\phi))=
\tilde{h}_{\tilde{A}(\tM,\tN,\phi)}(\tH_{(\tM,\tN,\phi)})\] 
corresponding via the correspondence in Lemma \ref{qf} to the identity in $h_{\tH_{(\tM,\tN,\phi)}}(\tH_{(\tM,\tN,\phi)})$.  By applying $R$ we get an element $R(\Theta_{(\tM,\tN,\phi)}) \in \mathcal{C}_{\tB} \times_{\mathcal{C}_{\tD}}^{h} \mathcal{C}_{\tC}(R(\tH_{(\tM,\tN,\phi)}),R(\tilde{A}(\tM,\tN,\phi)))$.  In order to apply Lemma  \ref{Keller} and conclude that $R$ is a dg-quasi-equivalence we need to show (1) the existence of quasi-isomorphisms $h_{\tH(\tM,\tN,\phi)} \to \tilde{h}_{\tilde{A}(\tM,\tN,\phi)}$ , (2) for every $\tE \in \mathcal{P}_{\tA}$ that $\tE$ and $\tH_{R(\tE)}$ are homotopy equivalent and (3) that $R(\tH_{(\tM,\tN,\phi)})$ and $(\tM,\tN,\phi) \in  \mathcal{P}_{\tB} \times_{\mathcal{P}_{\tD}}^{h} \mathcal{P}_{\tC}$ are homotopy equivalent for each $(\tM,\tN,\phi) \in  \mathcal{P}_{\tB} \times_{\mathcal{P}_{\tD}}^{h} \mathcal{P}_{\tC}$ via the map given in Lemma \ref{Keller}.  The first and second statement was already shown in Lemma \ref{injectivity}.  So, it only remains to show the third statement. 
We wish to show the composition 
\[\Lambda_{(\tM,\tN,\phi)}=\eta_{(\tM,\tN,\phi)} \circ R(\Theta_{(\tM,\tN,\phi)})
\] of $R(\Theta_{(\tM,\tN,\phi)})$ with $\eta_{(\tM,\tN,\phi)} \in \mathcal{C}_{\tB} \times_{\mathcal{C}_{\tD}}^{h} \mathcal{C}_{\tC}(R(\tilde{A}(\tM,\tN,\phi)),  (\tM,\tN,\phi))$ is a homotopy equivalence in $ \mathcal{P}_{\tB} \times_{\mathcal{P}_{\tD}}^{h} \mathcal{P}_{\tC}(R(\tH_{(\tM,\tN,\phi)}),  (\tM,\tN,\phi))$.  

The diagram
\begin{displaymath}
    \xymatrix{
        R(\tH_{(\tM,\tN,\phi)})  \ar[drr]_{\Lambda_{(\tM,\tN,\phi)}} \ar[rr]^{R(\Theta_{(\tM,\tN,\phi)})} & & R (\tilde{A}(\tM,\tN,\phi))  \ar[d]^{\eta_{(\tM,\tN,\phi)}}\\
                             &  & (\tM,\tN,\phi) }
\end{displaymath}
has as its zero component the diagram (\ref{zeroDiagOne}) applied to $(\tM^0,\tN^0,\phi^0)$.  Therefore, if we let $\Lambda_{(\tM,\tN,\phi)} = (\mu,\nu,\tau)$, we see that $\mu^0$ and $\nu^0$ are quasi-isomorphisms.  Therefore, by Proposition \ref{characterization} we see that $\mu$ and $\nu$ give isomorphisms in the homotopy categories $H^0( \mathcal{P}_\tB)$ and $H^0( \mathcal{P}_\tC)$.  Hence by Lemma \ref{getback}, we can conclude that the map $R(\tH_{(\tM,\tN,\phi)})  \to (\tM,\tN,\phi)$ is a homotopy equivalence.
\ \hfill $\Box$

\section{The Case of Complex Manifolds}\label{CCM}
We give here an example using the Dolbeault algebra, in the case $k=\mathbb{C}$.
Let $X$ be a compact complex manifold.   It was shown in \cite{B}, that there is an equivalence of triangulated categories between $\Ho\Pc_{\As(X)}$ and ${D}^b_{\mbox{coh}}(X)$. Even though it wasn't stated in \cite{B}, it is not too hard to derive the following stronger result. 
\begin{thm}\label{cateq} \cite{B} Let $X$ be a compact complex manifold and $\As(X)=(\A^\bullet(X), d, 0)=(\A^{0,\bullet}(X),\dbar, 0)$ the Dolbeault 
dga. Then there is a dq-quasi-equivalence of categories between $\Pc_{\As(X)}$ and $L_{pe}X$, where $L_{pe}X$ is the dg-enhancement of the category of perfect complexes of sheaves of $\mathcal{O}_{X}$ modules  as defined in \cite{T}.   
\end{thm}
\begin{rem} A holomorphic vector bundle $E$ on $X$ corresponds under the above equivalence to the pair $(\mE,\Eb) \in \mP_{\tA}$ where $\mE=H^{0}(X,E\otimes_{\mathcal{O}}\mathcal{C}^{\infty})$ is the space of global sections of the corresponding smooth vector bundle  thought of as an $\mA^0$ module and the connection $\mathbb{E}=\mathbb{E}^{1}$ is taken to be the obvious $\overline{\partial}$ connection on $\mE \otimes_{\mA^{0}} \mA$. 
\end{rem}
\begin{defn}
\begin{itemize}
Let $X$ be a complex manifold and $S \subset X$ any subset.   Define \item 
$\tA(S) = \varinjlim \tA(V)$
\item 
$
\stackrel{\circ}{\tA}(S)=\varinjlim\tA(V-S)
$
\end{itemize}
where in both cases the direct limit is taken over all open subspaces $V$ in the classical topology on $X$ which contain $S$.
\end{defn}

As a consequence of Theorem \ref{main} we have:
\begin{thm}
Let $X$ be a complex manifold.  Consider the following situations:

\begin{enumerate}[(a)]
\item $U_{1}$ and $U_{2}$ form an open cover of $X$.
\item $Z$ is a closed subspace of $X$ and $U=X-Z$

\end{enumerate}
Then the natural restriction functors are dg quasi-equivalences:
\begin{enumerate}[(a)]
\item \[
\Pc_{\As(X)}\to \Pc_{\As(\overline{U_1})}\times^h_{ \Pc_{\As(\overline{U_1}\cap \overline{U_2})} }\Pc_{\As(\overline{U_2})}
\]

\item
\[
\Pc_{\As(X)}\to \Pc_{\As(Z)}\times^h_{ \Pc_{\tD}}\Pc_{\As(U)}.\]
\end{enumerate}
where in (b) we take $\tD=\stackrel{\circ}{\tA}(Z)$
\end{thm}
{\bf Proof.}
We need to prove that the assumptions \ref{conditions} of Theorem \ref{main} hold. \begin{enumerate}[(a)]
\item
By definition of the direct limit, an element in $\As(\overline{U_1}\cap \overline{U_2})$ can be extended to some neighborhood of $\overline{U_1}\cap \overline{U_2}$ in $\overline{U_2}$ and thus it follows that any element of $\A(\overline{U_1}\cap \overline{U_2})$ can be extended to $\overline{U_2}$ because the Dolbeault sheaf is soft.    Therefore, the restriction $\tA(\overline{U_2}) \to \tA(\overline{U_1}\cap \overline{U_2})$ is surjective.  This is property (2).   Similarly, one can extend any pair of sections in $\mA^0(\overline{U_1})$ and  $\mA^0(\overline{U_2})$ which agree on the overlap to a unique element in $ \mA^0(X)$.  This is property (1).    Property (3) holds since the Dolbeault algebra on $X$ consists of the global sections of a vector bundle and therefore is a projective $\mA^{0}(X)$ module and hence a flat $\mA^{0}(X)$ module.    Moreover, the sections of a vector bundle satisfy property (4) since the sections of a vector bundle form a locally free sheaf over the smooth functions. 
\item
An element in $\stackrel{\circ}{\tA}(Z)$ can be extended to an element in $\tA(V-Z)$ where $V$ is an open neighborhood of $Z$.  By shrinking $V$ to another open set containing $Z$, one can extend some restriction of the element in $\tA(V-Z)$ to some element of $\tA(U)$.  This shows property (2).  The rest of the properties are routinely verified in a manner similar to (a).
\end{enumerate}
\ \hfill $\Box$
\begin{rem} For algebraic varieties over any field, there is a kind of descent similar to the second item where one uses the coherent sheaves on the formal neighborhood of $Z$ in place cohesive modules over $\tA(Z)$.  This can be found in \cite{BT}.  
\end{rem}

\section{A standard example}
It is well known that the failure of descent for derived categories is rectified by looking at dg-categories.  In this section we look at a standard example of this phenomenon in our context.

Let $U=\{z \in \mathbb{C} | |z| \leq 2 \}$ and $V=\{\xi \in \mathbb{C} | |\xi| \leq 2 \}$.  If we glue these along $W=\{z \in \mathbb{C} | \frac{1}{2} \leq |z| \leq 2 \}$ via $\xi=z^{-1}$ to form $\mathbb{P}^{1}_{\mathbb{C}}$ then the natural functor 
\[R: D^{b}(\mathbb{P}^{1}_{\mathbb{C}}) \to D^{b}(U) \times_{D^{b}(W)} D^{b}(V) 
\]
is not an equivalence of categories because a non-trivial element in  
\[ H^{1}(\mathbb{P}^{1}_{\mathbb{C}},\mathcal{O}(-2)) = \Ext^{1}_{\mathbb{P}^{1}_{\mathbb{C}}}(\mathcal{O}(2),\mathcal{O}) = \Hom_{D_{coh}^{b}(\mathbb{P}^{1}_{\mathbb{C}})}(\mathcal{O}(2),\mathcal{O}[1])
\]
is a morphism in the derived category which pulls back to a trivial morphism in the fiber product category $D^{b}(U) \times_{D^{b}(W)} D^{b}(V)$.
 Let $\tA, \tA(U), \tA(V), \tA(W)$ be the Dolbeaut dg-algebras of $\mathbb{P}^{1}_{\mathbb{C}}, U, V, W$ respectively.  Let the dg-category $\mQ$ be defined by \[\mQ= \mathcal{P}_{\tA(U)} \times_{\mathcal{P}_{\tA(W)}}^{h} \mathcal{P}_{\tA(V)}.\]  The information in the extension class is not lost via restriction in the dg context, in the sense that the natural map 
 \[R: \mathcal{P}_{\tA} \to \mQ
 \]
 preserves the information in the spaces of morphisms.   Let us see explicitly how this works in an example.  Consider $\tA(U), \tA(V)$ and $\tA(W)$ as objects of $\mathcal{P}_{\tA(U)}, \mathcal{P}_{\tA(V)}$ and $\mathcal{P}_{\tD}$ in the obvious way.  Notice that we have $\tg^{*}\tA(U)=\tA(W)=\tl^{*}\tA(V)$.  The line bundles $\mathcal{O}(j)$ on $\mathbb{P}^{1}$ correspond to the gluings given by the closed element $\phi \in \mathcal{P}^{0}_{\tA(W)}(\tg^{*}\tA(U),\tl^{*}\tA(V))$ which is  multiplication by $z^{-j} \in \mD^0$.  Let us define $\tA(j)$ by 
 \[\tA(j) = (\tA(U),\tA(V),z^{-j}) \in \mathcal{P}_{\tA(U)} \times_{\mathcal{P}_{\tA(W)}}^{h} \mathcal{P}_{\tA(V)}.
 \]
 We claim there is a natural isomorphism
\begin{equation}\label{weclaim}\Ext^{1}_{\mathbb{P}^{1}_{\mathbb{C}}}(\mathcal{O}(2),\mathcal{O}) \to H^{1}(\mQ(\tA(2),\tA)).
\end{equation}
In fact, we can exhibit a quasi-isomorphism 
\begin{equation}\label{exhibit}\check{C}^{\bullet}(\mathcal{H}om(\mathcal{O}(2),\mathcal{O})) \to \mQ^{\bullet}(\tA(2),\tA)
\end{equation}
between the \v{C}ech complex $\check{C}^{\bullet}(\mathcal{H}om(\mathcal{O}(2),\mathcal{O}))$ with respect to the cover by $U$ and $V$ and the complex $\mQ^{\bullet}(\tA(2),\tA)$.  The \v{C}ech complex $\check{C}^{\bullet}(\mathcal{H}om(\mathcal{O}(2),\mathcal{O}))$ with respect to the cover by $U$ and $V$ computes $\check{H}^{\bullet}(\mathcal{H}om(\mathcal{O}(2),\mathcal{O})) \cong \Ext^{\bullet}_{\mathbb{P}^{1}_{\mathbb{C}}}(\mathcal{O}(2),\mathcal{O})=\mathbb{C}$.  Notice that \[\check{C}^{1}(\mathcal{H}om(\mathcal{O}(2),\mathcal{O}))= \check{Z}^{1}(\mathcal{H}om(\mathcal{O}(2),\mathcal{O})) = \mathcal{O}(W).\]   Let $\tau \in \mathcal{O}(W)$ represent a generator of the one dimensional complex vector space 
\[\check{H}^{1}(\mathbb{P}^{1}_{\mathbb{C}},\mathcal{H}om(\mathcal{O}(2),\mathcal{O}))= \mathcal{O}(W)/\sim.\]  Choose a pair of smooth functions $(\mu_{U},\mu_{V}) \in C^{\infty}_{\mathbb{C}}(U) \times C^{\infty}_{\mathbb{C}}(V)$ so that 
\[z^{2}\mu_{U} - \mu_{V}= \tau.\]  
In degree one the map in (\ref{exhibit}) is defined by sending $\tau$ to the triple $(\overline{\partial}\mu_{U}, \overline{\partial} \mu_{V},0)$ thought of as a pair of degree one maps $\mA^{0}(U) \to\mA^{1}(U)$ and from $\mA^{0}(V) \to\mA^{1}(V)$  given by the product with $\overline{\partial}\mu_{U}$ and $\overline{\partial}\mu_{V}$ respectively along with the trivial homotopy between their restrictions to $W$.  It is easy to see that the corresponding class in $H^{1}(\mQ(\tA(2),\tA))$ is non-zero and does not depend on the choice of $\mu_{U}$ and $\mu_{V}$.  So we have produced a map as in (\ref{weclaim}) which is injective.  In order to see that it is surjective, consider a pair of closed degree one maps $\zeta_{U}: \mA^0(U) \to \mA^1(U), \zeta_{V}: \mA^0(V) \to \mA^1(V) $ and a homotopy $\gamma: \mA^0(W) \to \mA^1(W)$ between their restrictions $\zeta_{U}\otimes_{\mA^0(U)} \text{id}_{ \mA^0(W)}: \mA^0(W) \to  \mA^1(W)$ and 
$\zeta_{V} \otimes_{\mA^0(V)} \text{id}_{\mA^0(W)}:  \mA^0(W) \to \mA^1(W)$.  This triple is a general element in $Z^{1}(\mQ^{\bullet}(\tA(2),\tA))$.  We can think of these as $(0,1)$ forms $\zeta_{U}$ on $U$ and $\zeta_{V}$ on $V$ and a smooth function $\gamma$ on $W$ such that \[(z^{2}\zeta_{U}-\zeta_{V})|_{W}= \overline{\partial}\gamma.\]  The corresponding cohomology class in  $H^{1}(\mQ^{\bullet}(\tA(2),\tA))$ comes from the element 
\[z^{2} \rho_{U}-\rho_{V} -\gamma \in \mathcal{O}(W) = \check{Z}^{1}(\mathcal{H}om(\mathcal{O}(2),\mathcal{O})) \to \check{H}^{1}(\mathcal{H}om(\mathcal{O}(2),\mathcal{O}))
\] where we chose $\rho_{U}$ and $\rho_{V}$ such that $\overline{\partial}\rho_{U}= \zeta_{U}$ and $\overline{\partial}\rho_{V}= \zeta_{V}$. It is easy to see that the isomorphism (\ref{weclaim}) which we defined does not depend on the choice of $\tau$.  
\bibliographystyle{amsalpha}

\end{document}